\documentclass[11pt,reqno]{amsart}
\usepackage{amsfonts}
\usepackage{amssymb,latexsym}
\usepackage{cite}
\usepackage[height=190mm,width=143mm]{geometry}
\usepackage{lscape}

\theoremstyle{plain}
\newtheorem{theorem}{Theorem}

\theoremstyle{definition}

\theoremstyle{remark}

 \numberwithin{equation}{section}
\input{tcilatex}

\begin{document}

\begin{abstract}
In this article, we find bases for the spaces of modular forms $M_{3}(\Gamma
_{0}(40),\left( \frac{d}{\cdot }\right) )$ for $d=-4,-8,-20\text{ and }-40.$
We then derive formulas for the number of representations of a positive
integer by all the diagonal sextenary quadratic forms with coefficients $%
1,2,5$ and $10$.
\end{abstract}

\title{Representations by sextenary quadratic forms with coefficients $1$, $%
2 $, $5$ and $10$}
\author{B\"{u}lent K\"{o}kl\"{u}ce}
\address{The Institute for Computational and Experimental Research in
Mathematics\\
121 S Main St, Providence, 02903\\
Brown University\\
Rhode Island\\
USA}
\email{bulent\_kokluce@brown.edu}
\maketitle

\section{Introduction}

\bigskip Let $\mathbf{%
\mathbb{N}
,}$ $%
\mathbb{Z}
,$ $%
\mathbb{C}
$ and $%
\mathbb{H}
$ denote the set of positive integers, integers, complex numbers and the
upper half plane, respectively and let $\mathbf{%
\mathbb{N}
}_{\mathbf{0}}=\mathbf{%
\mathbb{N}
\cup \{0\}.}$Let $N\in \mathbf{%
\mathbb{N}
}$ and $\chi $ be a Dirichlet character of modulus dividing $N$ and $\Gamma
_{0}(N)$ the modular subgroup defined by

\begin{equation}
\Gamma _{0}(N)=\left\{ \left(
\begin{array}{cc}
a & b \\
c & d%
\end{array}%
\right) \in SL_{2}(%
\mathbb{Z}
):c\equiv 0(\text{mod }N)\right\} .  \label{1.3}
\end{equation}

Let $k\in
\mathbb{N}
.$ We write $M_{k}(\Gamma _{0}(N),\chi )$ to denote the space of modular
forms of weight $k$ with multiplier system $\chi $ for $\Gamma _{0}(N)$ and,
$E_{k}(\Gamma _{0}(N),\chi )$ and $S_{k}(\Gamma _{0}(N),\chi )$ to denote
the subspaces of Eisenstein forms and cusp forms of $M_{k}(\Gamma
_{0}(N),\chi ),$ respectively. It is known (see for example \cite[Theorem
2.1.7]{Miyake}{,} \cite[p.83]{Stein}) that

\begin{equation}
M_{k}(\Gamma _{0}(N),\chi )=E_{k}(\Gamma _{0}(N),\chi )\oplus S_{k}(\Gamma
_{0}(N),\chi ).  \label{e0}
\end{equation}

Let $a_{i},m\in \mathbf{%
\mathbb{N}
},r_{i},n\in \mathbf{%
\mathbb{N}
}_{\mathbf{0}}$ for all $1\leq i\leq m$ and let

\begin{equation}
N(a_{1}^{r_{1}},..,a_{m}^{r_{m}};n)  \label{1.1}
\end{equation}
denote the number of representations of $n$ by the quadratic form

\begin{equation}
\mathop{\displaystyle \sum }_{i=1}^{m}\mathop{\displaystyle \sum }\limits%
_{j=1}^{r_{i}}a_{i}x_{j}^{2}.  \label{1.1b}
\end{equation}%
Without loss of generality we may suppose that
\begin{equation}
a_{1}\leq ...\leq a_{m}\text{ }  \label{con}
\end{equation}%
For convenience, let us denote the sextenary quadratic form
\begin{eqnarray}
&&x_{1}^{2}+...+x_{l_{1}}^{2}+2(x_{l_{1}+1}^{2}+...+x_{l_{1}+l_{2}}^{2})+5(x_{l_{1}+l_{2}+1}^{2}+...+x_{l_{1}+l_{2}+l_{5}}^{2})
\label{1.2} \\
&&+10(x_{l_{1}+l_{2}+l_{5}+1}^{2}+...+x_{l_{1}+l_{2}+l_{5}+l_{10}}^{2})
\notag
\end{eqnarray}%
by $(l_{1},l_{2},l_{5},l_{10})$ and its representation number by $%
N(1^{l_{1}},2^{l_{2}},5^{l_{5}},10^{l_{10}};n).$

In this article, we find bases for the spaces of modular forms $M_{3}(\Gamma
_{0}(40),\left( \frac{d}{\cdot }\right) )$ for $d=-4,-8,-20\text{ and }-40.$
We then obtain some new theta function identities and use these results to
derive formulae $N(1^{i},2^{j},5^{k},10^{l};n)$ for each of the sextenary
quadratic forms $(l_{1},l_{2},l_{5},l_{10})$ with coefficients $1,2,5$ or $10
$. There is a total of $84$ sextenary quadratic forms with these
coefficients(see Table \ref{Tab:Tb1}). To the best of our knowledge, among
these $84$ formulae, except $7$ of them which appeared in the literature,
all the rest are new. Formulae for $N(1^{6},2^{0},5^{0},10^{0};n),$ $%
N(1^{2},2^{4},5^{0},10^{0};n)$ and $N(1^{4},2^{2},5^{0},10^{0};n)$ are given
in \cite{Alaca1}. Aygin \cite{Aygin} found formulae for each diagonal
sextenary quadratic form with coefficients $1,2,3,$ and $6$. In the current
study and in his work, since the cases with $l_{1}+l_{2}=6$ are common, we
compared the results and see that they are all consistent.

For some other works on representation numbers of sextenary quadratic forms
with other coefficients, see \cite{Alaca2,Alaca3,Alaca4,Berkovich,Xia}. See
\cite{Grosswald,Williams} for a classical history of this research and \cite%
{Alaca5,Alaca6,Cooper,Kokluce,Kokluce2,Milne,Ono,Ramakrishnan} for
contemporary accounts of the subject.

\begin{table}[th]
\caption{List of sextenary quadratic forms $(l_{1},l_{2},l_{5},l_{10})$}
\label{Tab:Tb1}\centering
\begin{tabular}{|cccc|}
\hline
$(l_{1},l_{2},l_{5},l_{10})$ & $(l_{1},l_{2},l_{5},l_{10})$ & $%
(l_{1},l_{2},l_{5},l_{10})$ & $(l_{1},l_{2},l_{5},l_{10})$ \\ \hline\hline
$(6,0,0,0)$ $(4,2,0,0)$ & $(5,0,1,0)$ $(4,1,0,1)$ & $(5,1,0,0)$ $(4,0,1,1)$
& $(5,0,0,1)$ $(4,1,1,0)$ \\
$(4,0,2,0)$ $(4,0,0,2)$ & $(3,2,1,0)$ $(3,0,3,0)$ & $(3,3,0,0)$ $(3,1,2,0)$
& $(3,2,0,1)$ $(3,0,2,1)$ \\
$(3,1,1,1)$ $(2,4,0,0)$ & $(3,0,1,2)$ $(2,3,0,1)$ & $(3,1,0,2)$ $(2,2,1,1)$
& $(3,0,0,3)$ $(2,3,1,0)$ \\
$(2,2,2,0)$ $(2,2,0,2)$ & $(2,1,2,1)$ $(2,1,0,3)$ & $(2,0,3,1)$ $(2,0,1,3)$
& $(2,1,3,0)$ $(2,1,1,2)$ \\
$(2,0,4,0)$ $(2,0,2,2)$ & $(1,4,1,0)$ $(1,2,3,0)$ & $(1,5,0,0)$ $(1,3,2,0)$
& $(1,4,0,1)$ $(1,2,2,1)$ \\
$(2,0,0,4)$ $(1,3,1,1)$ & $(1,2,1,2)$ $(1,0,5,0)$ & $(1,3,0,2)$ $(1,1,4,0)$
& $(1,2,0,3)$ $(1,0,4,1)$ \\
$(1,1,3,1)$ $(1,1,1,3)$ & $(1,0,3,2)$ $(1,0,1,4)$ & $(1,1,2,2)$ $(1,1,0,4)$
& $(1,0,2,3)$ $(1,0,0,5)$ \\
$(0,6,0,0)$ $(0,4,2,0)$ & $(0,5,0,1)$ $(0,3,2,1)$ & $(0,4,1,1)$ $(0,2,3,1)$
& $(0,5,1,0)$ $(0,3,3,0)$ \\
$(0,4,0,2)$ $(0,2,4,0)$ & $(0,3,0,3)$ $(0,1,4,1)$ & $(0,2,1,3)$ $(0,0,5,1)$
& $(0,3,1,2)$ $(0,1,5,0)$ \\
$(0,2,2,2)$ $(0,2,0,4)$ & $(0,1,2,3)$ $(0,1,0,5)$ & $(0,0,3,3)$ $(0,0,1,5)$
& $(0,1,3,2)$ $(0,1,1,4)$ \\
$(0,0,6,0)$ $(0,0,4,2)$ &  &  &  \\
$(0,0,2,4)$ $(0,0,0,6)$ &  &  &  \\ \hline
\end{tabular}%
\end{table}

\section{Preliminary Results}

\bigskip The Dedekind eta function $\eta (z)$ is the holomorphic function
defined on the upper half plane $%
\mathbb{H}
$ $=\left\{ z\in
\mathbb{C}
:\func{Im}(z)>0\right\} $ by the product formula

\begin{equation}
\eta (z)=e^{\pi iz/12}\mathop{\displaystyle \prod }\limits_{n=1}^{\infty
}(1-e^{2\pi inz}).  \label{kk}
\end{equation}
Through the remainder of the paper we take $q=q(z):=e^{2\pi iz}$ with $z$ in
$%
\mathbb{H}
$ and so by (\ref{kk}) we have
\begin{equation}
\eta (z)=q^{1/24}\mathop{\displaystyle \prod }\limits_{n=1}^{\infty
}(1-q^{n}).  \label{e1}
\end{equation}
An eta quotient is defined to be a finite product of the form
\begin{equation}
f(z)=\mathop{\displaystyle \prod }\limits_{\delta }\eta ^{r_{\delta
}}(\delta z),  \label{e2}
\end{equation}
where $\delta $ runs through a finite set of positive integers and the
exponents $r_{\delta }$ are nonzero integers. By taking $N$ to be the least
common multiple of $\delta ^{\prime }s$ we can write the eta quotient (\ref%
{e2}) as
\begin{equation}
f(z)=\mathop{\displaystyle \prod }\limits_{\delta \mid N}\eta ^{r_{\delta
}}(\delta z).  \label{e3}
\end{equation}

For $q\in
\mathbb{C}
$ with $\left\vert q\right\vert <1$, the Ramanujan theta function $\varphi
(q)$ is defined by

\begin{equation}
\varphi (q)=\mathop{\displaystyle \sum }\limits_{n=-\infty }^{\infty
}q^{n^{2}}.  \label{eta}
\end{equation}%
We have
\begin{equation}
\mathop{\displaystyle \sum }\limits_{n=0}^{\infty
}N(1^{l_{1}},2^{l_{2}},5^{l_{5}},10^{l_{10}};n)q^{n}=\mathop{\displaystyle
\prod }\limits_{d\mid 10}\varphi ^{l_{d}}(q^{d}).  \label{tet1}
\end{equation}%
It is well known by Jacobi's triple product identity \cite[p.10]{Berndt} that

\begin{equation}
\varphi (q)=\frac{\eta ^{5}(2z)}{\eta ^{2}(z)\eta ^{2}(4z)}.  \label{tet3}
\end{equation}
So, we can rewrite the generating function (\ref{tet1}) in terms of
eta quotients as follows:

\begin{equation*}
\mathop{\displaystyle \sum }\limits_{n=0}^{\infty
}N(1^{l_{1}},2^{l_{2}},5^{l_{5}},10^{l_{10}};n)q^{n}=\mathop{\displaystyle
\prod }\limits_{d\mid 10}\left( \frac{\eta ^{5}(2dz)}{\eta ^{2}(dz)\eta
^{2}(4dz)}\right) ^{l_{d}}.
\end{equation*}

Let $\chi $ and $\psi $ be Dirichlet characters. For $n\in
\mathbb{N}
$ we define $\sigma _{(k,\chi ,\psi )}$ by

\begin{equation}
\sigma _{(k,\chi ,\psi )}(n):=\sum_{1\leq d\mid n}\psi (d)\chi (n/d)d^{k}.
\label{sig}
\end{equation}

\bigskip If $n\notin
\mathbb{N}
$ we set $\sigma _{(k,\chi ,\psi )}(n)=0.$ For each quadratic discriminant $%
t $, we put $\chi _{t}(n)=\left( \frac{t}{n}\right) ,$ where $\left( \frac{t%
}{n}\right) $ is the Kronecker symbol. Now we define the following five
Dirichlet characters:

\begin{equation}
\left\{
\begin{tabular}{lll}
$\chi _{-4}(n)=\left( \frac{-4}{n}\right) ,$ & $\chi _{-8}(n)=\left( \frac{-8%
}{n}\right) ,$ & $\chi _{5}(n)=\left( \frac{5}{n}\right) ,$ \\
&  &  \\
$\chi _{-20}(n)=\left( \frac{-20}{n}\right) ,$ & $\chi _{-40}(n)=\left(
\frac{-40}{n}\right) .$ &
\end{tabular}%
\right.  \label{char}
\end{equation}

Let $\chi $ and $\psi $ be primitive Dirichlet characters with conductors $%
L,R\in
\mathbb{N}
$ respectively such that $\chi (-1)\psi (-1)=-1$. The weight $3$ Eisenstein
series are defined by

\begin{equation}
E_{3,\chi ,\psi }(q)=c_{0}+\mathop{\displaystyle \sum }\limits_{n\geq
1}\sigma _{(2,\chi ,\psi )}(n)q^{n}  \label{eis}
\end{equation}%
where
\begin{equation}
c_{0}=\left\{
\begin{array}{c}
-\frac{B_{3},_{\chi }}{6} \\
0%
\end{array}%
\right.
\begin{array}{c}
\text{if }L=1 \\
\text{if }L>1%
\end{array}
\label{eis2}
\end{equation}%
where the generalized Bernoulli numbers $B_{3,\chi }$ attached to $\chi $
are defined by the following equation:
\begin{equation}
B_{3,\chi }=6[x^{3}]\mathop{\displaystyle \sum }\limits_{a=1}^{L}\frac{\chi
(a)xe^{ax}}{e^{Lx}-1}.  \label{eis3}
\end{equation}

From (\ref{eis3}) we compute
\begin{equation}
B_{3,\chi _{-4}}=\frac{3}{2},\text{ }B_{3,\chi _{-20}}=90,B_{3,\chi _{-8}}=9,%
\text{ }B_{3,\chi _{-40}}=474.
\end{equation}%
From \cite[Theorem 5.8, p.88]{Stein} we know that except when $k=2$ and $%
\chi =\psi =1$, the Eisenstein series $E_{k,\chi ,\psi }$ defines an element
of $M_{k}(\Gamma _{0}(RLt),\chi \psi )$ where $R$ and $L$ are respectively
the conductors of $\chi $and $\psi .$ Further from \cite[Theorem 5.9, p.88]%
{Stein}, the Eisenstein series obtained from here with $RLt\mid N$ and $\chi
\psi =\varepsilon $ form a basis for the Eisenstein subspace $E_{k}(\Gamma
_{0}(N),\varepsilon )$. From (\ref{eis}) we can write weight $3$ level $40$
Eisenstein series as follows:

\begin{equation}
E_{3,\chi _{1},\chi _{-4}}(q):=-\frac{1}{4}+\mathop{\displaystyle \sum }%
\limits_{n=1}^{\infty }\sigma _{(2,\chi _{1},\chi _{-4})}(n)q^{n},\text{ }%
E_{3,\chi _{-4},\chi _{1}}(q):=\mathop{\displaystyle \sum }%
\limits_{n=1}^{\infty }\sigma _{(2,\chi _{-4},\chi _{1})}(n)q^{n},
\label{Eis2}
\end{equation}

\begin{equation}
E_{3,\chi _{1},\chi _{-20}}(q):=-15+\mathop{\displaystyle \sum }%
\limits_{n=1}^{\infty }\sigma _{(2,\chi _{1},\chi _{-20})}(n)q^{n},\text{ }%
E_{3,\chi _{-20},\chi _{1}}(q):=\mathop{\displaystyle \sum }%
\limits_{n=1}^{\infty }\sigma _{(2,\chi _{-20},\chi _{1})}(n)q^{n},
\label{Eis3}
\end{equation}%
\begin{equation}
E_{3,\chi _{-4},\chi _{5}}(q):=\mathop{\displaystyle \sum }%
\limits_{n=1}^{\infty }\sigma _{(2,\chi _{-4},\chi _{5})}(n)q^{n},\text{ }%
E_{3,\chi _{5},\chi _{-4}}(q):=\mathop{\displaystyle \sum }%
\limits_{n=1}^{\infty }\sigma _{(2,\chi _{5},\chi _{-4})}(n)q^{n},
\label{Eis4}
\end{equation}

\begin{equation}
E_{3,\chi _{1},\chi _{-8}}(q):=-\frac{3}{2}+\mathop{\displaystyle \sum }%
\limits_{n=1}^{\infty }\sigma _{(2,\chi _{1},\chi _{-8})}(n)q^{n},\text{ }%
E_{3,\chi _{-8},\chi _{1}}(q):=\mathop{\displaystyle \sum }%
\limits_{n=1}^{\infty }\sigma _{(2,\chi _{-8},\chi _{1})}(n)q^{n},
\label{Eis5}
\end{equation}

\begin{equation}
E_{3,\chi _{1},\chi _{-40}}(q):=-79+\mathop{\displaystyle \sum }%
\limits_{n=1}^{\infty }\sigma _{(2,\chi _{1},\chi _{-40})}(n)q^{n},\text{ }%
E_{3,\chi _{-40},\chi _{1}}(q):=\mathop{\displaystyle \sum }%
\limits_{n=1}^{\infty }\sigma _{(2,\chi _{-40},\chi _{1})}(n)q^{n},
\label{Eis6}
\end{equation}%
\begin{equation}
E_{3,\chi _{-8},\chi _{5}}(q):=\mathop{\displaystyle \sum }%
\limits_{n=1}^{\infty }\sigma _{(2,\chi _{-8},\chi _{5})}(n)q^{n},\text{ }%
E_{3,\chi _{5},\chi _{-8}}(q):=\mathop{\displaystyle \sum }%
\limits_{n=1}^{\infty }\sigma _{(2,\chi _{5},\chi _{-8})}(n)q^{n}.
\label{Eis7}
\end{equation}%
M. Newman \cite{Newman,Newman2} systematically used the eta functions to
construct modular forms for $\Gamma _{0}(N)$ and then to check whether a
function $f(z)$ was a modular form for $\Gamma _{0}(N)$. The order of
vanishing of an eta function at the cusps of $\Gamma _{0}(N)$, was
determined by G. Ligozat \cite{Ligozat}. The following theorem which is
referred as Ligozat's Criteria (see \cite[Theorem 5.7, p.99]{Kilford},\cite[%
Corollary 2.3, p.37]{Kohler}) is used to determine if an eta quotient is in $%
M_{k}(\Gamma _{0}(N),\chi )$.

\begin{theorem}
\label{thm:1} Let $N$ be a positive integer and let $f(z)=%
\mathop{\displaystyle \prod }\nolimits_{1\leq \delta \mid N}\eta ^{r_{\delta
}}(\delta z)$ be an eta quotient. Let $s=\mathop{\displaystyle \prod }%
\nolimits_{1\leq \delta \mid N}\eta ^{\left\vert r_{\delta }\right\vert }.$
Suppose that the following conditions hold:

\begin{enumerate}
\item[(i)] $\mathop{\displaystyle \sum }\nolimits_{1\leq \delta \mid
N}\delta .r_{\delta }\equiv 0(\func{mod}24),$

\item[(ii)] $\mathop{\displaystyle \sum }\nolimits_{1\leq \delta \mid N}%
\dfrac{N}{\delta }.r_{\delta }\equiv 0(\func{mod}24),$

\item[(iii)] the weight $k=\dfrac{1}{2}\mathop{\displaystyle \sum }%
\nolimits_{1\leq \delta \mid N}r_{\delta }$ is an even integer,

\item[(iv)] for each $d\mid N,\mathop{\displaystyle \sum }\nolimits_{1\leq
\delta \mid N}\dfrac{\gcd (d,\delta )^{2}.r_{\delta }}{\delta }\geq 0.$
\end{enumerate}
\end{theorem}

Then, $f(z)$ is in $M_{k}(\Gamma _{0}(N),\chi ),$ where the character $\chi $
is defined by $\chi (m)=(\frac{(-1)^{k}s}{m}).$ In addition to the above
conditions, if all the inequalities in (iv) hold strictly, then $f(z)$ is in
$S_{k}(\Gamma _{0}(N),\chi ).$

Using the dimension formulas given in \cite[p.98]{Stein} we obtain the
dimensions of the spaces of Eisenstein forms and cusp forms for each
character $\chi .$ We see that
\begin{equation}
\text{dim}(E_{3}(\Gamma _{0}(40),\chi _{-4})=8,\text{dim}(S_{3}(\Gamma
_{0}(40),\chi _{-4}))=8.  \label{dim1}
\end{equation}%
\begin{equation}
\text{dim}(E_{3}(\Gamma _{0}(40),\chi _{-20}))=8,\text{ dim}(S_{3}(\Gamma
_{0}(40),\chi _{-20}))=8  \label{dim2}
\end{equation}%
\begin{equation}
\text{dim}(E_{3}(\Gamma _{0}(40),\chi _{-8}))=4,\text{dim}(S_{3}(\Gamma
_{0}(40),\chi _{-8}))=10,  \label{dim3}
\end{equation}%
and
\begin{equation}
\text{dim}(E_{3}(\Gamma _{0}(40),\chi _{-40}))=4,\text{ dim}(S_{3}(\Gamma
_{0}(40),\chi _{-40}))=10.  \label{dim4}
\end{equation}

\section{Bases for\ the Spaces of Modular Forms}

\subsection{Basis for $M_{3}(\Gamma _{0}(40),\protect\chi _{-4})$}

We define the following $8$ eta quotients:

\begin{equation}
A_{1}(q)=\frac{\eta ^{4}(z)\eta ^{2}(4z)\eta ^{2}(20z)}{\eta ^{2}(2z)},
\label{a1}
\end{equation}

\begin{equation}
A_{2}(q)=\frac{\eta ^{3}(z)\eta (2z)\eta (5z)\eta (8z)\eta ^{3}(20z)}{\eta
(4z)\eta (10z)\eta (40z)},  \label{a2}
\end{equation}%
\begin{equation}
A_{3}(q)=\frac{\eta ^{3}(z)\eta ^{2}(4z)\eta (5z)\eta (40z)}{\eta (8z)},
\label{a3}
\end{equation}

\begin{equation}
A_{4}(q)=\frac{\eta ^{2}(z)\eta (4z)\eta (8z)\eta ^{4}(10z)\eta (20z)}{\eta
^{2}(5z)\eta (40z)},  \label{a4}
\end{equation}

\begin{equation}
A_{5}(q)=\frac{\eta ^{2}(z)\eta ^{2}(5z)\eta ^{4}(20z)}{\eta ^{2}(10z)},
\label{a5}
\end{equation}

\begin{equation}
A_{6}(q)=\frac{\eta (z)\eta (2z)\eta ^{2}(8z)\eta (10z)\eta ^{4}(20z)}{\eta
(5z)\eta ^{2}(40z)},  \label{a6}
\end{equation}

\begin{equation}
A_{7}(q)=\frac{\eta (z)\eta (8z)\eta ^{2}(10z)\eta ^{3}(40z)}{\eta (5z)},
\label{a7}
\end{equation}

\begin{equation}
A_{8}(q)=\frac{\eta (2z)\eta ^{2}(8z)\eta ^{3}(20z)\eta ^{2}(40z)}{\eta
(4z)\eta (10z)},  \label{a8}
\end{equation}%
and the integers $a_{j}(n)$ $(n\in
\mathbb{N}
)$ by
\begin{equation}
A_{j}(q)=\mathop{\displaystyle \sum }\limits_{n=1}^{\infty }a_{j}(n)q^{n},%
\text{ }1\leq j\leq 8.  \label{a9}
\end{equation}

\begin{theorem}
\label{thm:2}$\{E_{3,\chi _{1},\chi _{-4}}(q^{t}),$ $E_{3,\chi _{-4,}\chi
_{1}}(q^{t}):t=1,2,5,10\}\cup \{A_{j}(q):j=1,2,...,8\}$ form a basis of $%
M_{3}(\Gamma _{0}(40),\chi _{-4}).$
\end{theorem}

\begin{proof}
From \cite[p.98]{Stein} we see that the\ dimension of $E_{3}(\Gamma
_{0}(40),\chi _{-4})$ is $8$. Appealing to \cite[Theorem 5.9, p.88]{Stein}
with $\epsilon =\chi _{-4}$ and $\psi ,$ $\chi $ $\in \{\chi _{1},$ $\chi
_{-4}\}$, we see that $\{E_{3,\chi _{1},\chi _{-4}}(q^{t}),E_{3,\chi
_{-4},\chi _{1}}(q^{t}):t=1,2,5,10\}$ is a basis for $E_{3}(\Gamma
_{0}(40),\chi _{-4}).$ From (\ref{e1})-(\ref{e3}), (\ref{a1})-(\ref{a8}),
and Theorem \ref{thm:1} we see that all the eta quotients $A_{j}(q)(1\leq
j\leq 8)$ are contained in $S_{3}(\Gamma _{0}(40),\chi _{-4}).$ The eta
quotients $A_{j}(q)(1\leq j\leq 10)$ are linearly independent over $%
\mathbb{C}
$. Since dim($S_{3}(\Gamma _{0}(40),\chi _{-4}))=8$ from\ (\ref{dim1})$,$ $%
\{A_{j}(q):j=1,2,...,8\}$ constitute a basis for $S_{3}(\Gamma _{0}(40),\chi
_{-4}).$

Thus $\{E_{3,\chi _{1},\chi _{-4}}(q^{t}),$ $E_{3,\chi _{-4},\chi
_{1}}(q^{t}):t=1,2,5,10\}\ $together with $\{A_{j}(q):j=1,2,...,8\}$ form a
basis for$M_{3}(\Gamma _{0}(40),\chi _{-4}).$
\end{proof}

\subsection{\protect\bigskip Basis for $M_{3}(\Gamma _{0}(40),\protect\chi %
_{-20})$}

We define the following $8$ eta quotients:

\begin{equation}
B_{1}(q)=\frac{\eta ^{4}(z)\eta ^{2}(2z)\eta (20z)}{\eta (4z)},  \label{b1}
\end{equation}

\begin{equation}
B_{2}(q)=\frac{\eta ^{4}(z)\eta ^{4}(20z)}{\eta (2z)\eta (10z)},  \label{b2}
\end{equation}

\begin{equation}
B_{3}(q)=\frac{\eta ^{3}(z)\eta (5z)\eta (8z)\eta (10z)\eta ^{2}(20z)}{\eta
(2z)\eta (40z)},  \label{b3}
\end{equation}

\begin{equation}
B_{4}(q)=\frac{\eta ^{2}(z)\eta (2z)\eta (8z)\eta ^{3}(10z)\eta ^{3}(20z)}{%
\eta (4z)\eta ^{2}(5z)\eta (40z)},  \label{b4}
\end{equation}

\begin{equation}
B_{5}(q)=\frac{\eta ^{2}(z)\eta ^{2}(4z)\eta ^{4}(10z)\eta (40z)}{\eta
^{2}(5z)\eta (8z)},  \label{b5}
\end{equation}

\begin{equation}
B_{6}(q)=\frac{\eta (z)\eta (2z)\eta (4z)\eta (10z)\eta ^{3}(20z)}{\eta (5z)}%
,  \label{b6}
\end{equation}

\begin{equation}
B_{7}(q)=\frac{\eta (z)\eta ^{2}(10z)\eta ^{4}(20z)}{\eta (5z)},  \label{b7}
\end{equation}

\begin{equation}
B_{8}(q)=\frac{\eta (2z)\eta ^{2}(20z)\eta ^{4}(40z)}{\eta (10z)},
\label{b8}
\end{equation}%
and the integers $b_{j}(n)$ $(n\in
\mathbb{N}
)$ by

\begin{equation}
B_{j}(q)=\mathop{\displaystyle \sum }\limits_{n=1}^{\infty }b_{j}(n)q^{n},%
\text{ }1\leq j\leq 8.  \label{b9}
\end{equation}

\begin{theorem}
\label{thm:3}$\{E_{3,\chi _{1},\chi _{-20}}(q^{t}),$ $E_{3,\chi _{-20},\chi
_{1}}(q^{t}),$ $E_{3,\chi _{-4},\chi _{5}}(q^{t}),$ $E_{3,\chi _{5},\chi
_{-4}}(q^{t}):t=1,2\}\cup \{B_{j}(q):j=1,2,...,8\}$ form a basis of $%
M_{3}(\Gamma _{0}(40),\chi _{-20}).$
\end{theorem}

\begin{proof}
Appealing to \cite[Theorem 5.9, p.88]{Stein} with $\epsilon =\chi _{-20}$
and $\psi ,$ $\chi $ $\in \{\chi _{1},$ $\chi _{-20},$ $\chi _{-4},$ $\chi
_{5}\}$, we see that $\{E_{3,\chi _{1},\chi _{-20}}(q^{t}),E_{3,\chi
_{-20},\chi _{1}}(q^{t}),E_{3,\chi _{-4},\chi _{5}}(q^{t}),E_{3,\chi
_{5},\chi _{-4}}(q^{t}):t=1,2\}$ is a basis for $E_{3}(\Gamma _{0}(40),\chi
_{-20}).$ From (\ref{e1})-(\ref{e3}), (\ref{b1})-(\ref{b8}) and Theorem \ref%
{thm:1} we see that all the eta quotients $B_{j}(q)(1\leq j\leq 8)$ are
contained in $S_{3}(\Gamma _{0}(40),\chi _{-20}).$ The eta quotients $%
B_{j}(q)(1\leq j\leq 8)$ are linearly independent over $%
\mathbb{C}
$. Since dim($S_{3}(\Gamma _{0}(40),\chi _{-20}))=8$ from\ (\ref{dim2}), we
see that $\{B_{j}(q):j=1,2,...,8\}$ form a basis of $S_{3}(\Gamma
_{0}(40),\chi _{-20})\ $and thus $\{E_{3,\chi _{1},\chi
_{-20}}(q^{t}),E_{3,\chi _{-20},\chi _{1}}(q^{t}),E_{3,\chi _{-4},\chi
_{5}}(q^{t}),E_{3,\chi _{5},\chi _{-4}}(q^{t}):t=1,2\}\ $together with $%
\{B_{j}(q):j=1,2,...,8\}$ form a basis of $M_{3}(\Gamma _{0}(40),\chi
_{-20}).$
\end{proof}

\subsection{\protect\bigskip Basis for $M_{3}(\Gamma _{0}(40),\protect\chi %
_{-8})$}

We define the following $10$ eta quotients:

\begin{equation}
C_{1}(q)=\frac{\eta ^{4}(z)\eta (2z)\eta ^{4}(10z)\eta ^{2}(40z)}{\eta
(4z)\eta ^{2}(5z)\eta ^{2}(20z)},  \label{c1}
\end{equation}

\begin{equation}
C_{2}(q)=\frac{\eta ^{4}(z)\eta ^{2}(5z)\eta (8z)\eta (40z)}{\eta (4z)\eta
(10z)},  \label{c2}
\end{equation}

\begin{equation}
C_{3}(q)=\frac{\eta ^{4}(z)\eta ^{3}(10z)\eta (20z)\eta ^{2}(40z)}{\eta
^{2}(2z)\eta ^{2}(5z)},  \label{c3}
\end{equation}

\begin{equation}
C_{4}(q)=\frac{\eta ^{3}(z)\eta ^{3}(8z)\eta ^{3}(10z)\eta (20z))}{\eta
^{2}(4z)\eta (5z)\eta (40z)},  \label{c4}
\end{equation}

\begin{equation}
C_{5}(q)=\frac{\eta ^{3}(z)\eta ^{3}(5z)\eta (20z)\eta ^{2}(40z)}{\eta
(2z)\eta ^{2}(10z)},  \label{c5}
\end{equation}

\begin{equation}
C_{6}(q)=\frac{\eta ^{3}(z)\eta (4z)\eta ^{3}(5z)\eta ^{2}(8z)}{\eta
^{2}(2z)\eta (10z)},  \label{c6}
\end{equation}

\begin{equation}
C_{7}(q)=\frac{\eta ^{2}(z)\eta (2z)\eta ^{3}(8z)\eta ^{3}(40z)}{\eta
^{2}(4z)\eta (20z)},  \label{c7}
\end{equation}

\begin{equation}
C_{8}(q)=\frac{\eta ^{2}(z)\eta (4z)\eta ^{4}(20z)\eta ^{2}(40z)}{\eta
(2z)\eta ^{2}(10z)},  \label{c8}
\end{equation}

\begin{equation}
C_{9}(q)=\frac{\eta (z)\eta (5z)\eta ^{2}(8z)\eta ^{4}(40z)}{\eta (4z)\eta
(10z)},  \label{c9}
\end{equation}

\begin{equation}
C_{10}(q)=\frac{\eta ^{2}(2z)\eta (4z)\eta ^{2}(8z)\eta ^{3}(10z)}{\eta
(z)\eta (5z)},  \label{c10}
\end{equation}%
and the integers $c_{j}(n)$ $(n\in
\mathbb{N}
)$ by

\begin{equation}
C_{j}(q)=\mathop{\displaystyle \sum }\limits_{n=1}^{\infty }c_{j}(n)q^{n},%
\text{ }1\leq j\leq 10.  \label{c11}
\end{equation}

\begin{theorem}
\label{thm:4}$\{E_{3,\chi _{_{1}},\chi _{_{-8}}}(q^{t}),E_{3,\chi _{-8},\chi
_{_{1}}}(q^{t}):t=1,5\}\cup \{C_{j}(q):j=1,2,...,10\}$ form a basis of $%
M_{3}(\Gamma _{0}(40),\chi _{-8}).$
\end{theorem}

\begin{proof}
Appealing to \cite[Theorem 5.9, p.88]{Stein} with $\epsilon =\chi _{-8},$
and $\psi ,$ $\chi \in \{\chi _{1},\chi _{-8}\}$, we see that $\{E_{3,\chi
_{_{1}},\chi _{_{-8}}}(q^{t}),E_{3,\chi _{-8},\chi _{_{1}}}(q^{t}):t=1,5\}$
is a basis for $E_{3}(\Gamma _{0}(40),\chi _{-8}).$ From (\ref{e1})-(\ref{e3}%
), (\ref{c1})-(\ref{c10}) and Theorem \ref{thm:1} we see that all the eta
quotients $C_{j}(q)(1\leq j\leq 10)$ are contained in $S_{3}(\Gamma
_{0}(40),\chi _{-8}).$ Since they are linearly independent and dim($%
S_{3}(\Gamma _{0}(40),\chi _{-8}))=10$ from (\ref{dim3})$,$ $%
\{C_{j}(q):j=1,2,...,10\}$ form a basis of $S_{3}(\Gamma _{0}(40),\chi
_{-8}).$

Thus, $\{E_{3,\chi _{_{1}},\chi _{_{-8}}}(q^{t}),E_{3,\chi _{-8},\chi
_{_{1}}}(q^{t}):t=1,5\}\ $together with $\{C_{j}(q):j=1,2,...,10\}$ form a
basis of $M_{3}(\Gamma _{0}(40),\chi _{-8}).$
\end{proof}

\subsection{\protect\bigskip Basis for $M_{3}(\Gamma _{0}(40),\protect\chi %
_{-40})$}

We define the following $10$ eta quotients:

\begin{equation}
D_{1}(q)=\frac{\eta ^{4}(z)\eta ^{2}(8z)\eta ^{3}(10z)\eta ^{2}(20z)}{\eta
^{2}(2z)\eta (4z)\eta ^{2}(5z)},  \label{d1}
\end{equation}

\begin{equation}
D_{2}(q)=\frac{\eta ^{3}(z)\eta (2z)\eta (8z)\eta ^{2}(10z)\eta (40z)}{\eta
(5z)\eta (20z)},  \label{d2}
\end{equation}

\begin{equation}
D_{3}(q)=\frac{\eta ^{3}(z)\eta ^{3}(5z)\eta ^{2}(8z)\eta ^{2}(20z)}{\eta
(2z)\eta (4z)\eta ^{2}(10z)},  \label{d3}
\end{equation}

\begin{equation}
D_{4}(q)=\frac{\eta ^{2}(z)\eta ^{2}(2z)\eta ^{3}(10z)\eta ^{2}(40z)}{\eta
(4z)\eta ^{2}(20z)},  \label{d4}
\end{equation}

\begin{equation}
D_{5}(q)=\frac{\eta ^{2}(z)\eta (2z)\eta ^{4}(5z)\eta (8z)\eta (40z)}{\eta
(4z)\eta ^{2}(10z)},  \label{d5}
\end{equation}

\begin{equation}
D_{6}(q)=\frac{\eta ^{2}(z)\eta (2z)\eta ^{2}(8z)\eta ^{3}(20z)}{\eta
^{2}(4z)},  \label{d6}
\end{equation}

\begin{equation}
D_{7}(q)=\frac{\eta ^{2}(z)\eta ^{2}(4z)\eta ^{2}(8z)\eta ^{2}(10z)}{\eta
(2z)\eta (20z)},  \label{d7}
\end{equation}

\begin{equation}
D_{8}(q)=\frac{\eta ^{2}(z)\eta (4z)\eta ^{2}(8z)\eta ^{3}(10z)}{\eta
^{2}(2z)},  \label{d8}
\end{equation}

\begin{equation}
D_{9}(q)=\frac{\eta (z)\eta (5z)\eta (8z)\eta ^{4}(20z)\eta (40z)}{\eta
(4z)\eta (10z)},  \label{d9}
\end{equation}

\begin{equation}
D_{10}(q)=\frac{\eta ^{2}(5)\eta ^{2}(8z)\eta (10z)\eta ^{2}(20z)}{\eta (4z)}%
,  \label{d10}
\end{equation}%
and the integers $d_{j}(n)$ $(n\in
\mathbb{N}
)$ by
\begin{equation}
D_{j}(q)=\mathop{\displaystyle \sum }\limits_{n=1}^{\infty }d_{j}(n)q^{n},%
\text{ }1\leq j\leq 10.  \label{d11}
\end{equation}

\begin{theorem}
\label{thm:5}$\{E_{3,\chi _{1},\chi _{-40}}(q),$ $E_{3,\chi _{-40},\chi
_{1}}(q),$ $E_{3,\chi _{-8},\chi _{5}}(q),$ $E_{3,\chi _{5},\chi
_{-8}}(q)\}\cup \{D_{j}(q):j=1,2,...,10\}$ form a basis of $M_{3}(\Gamma
_{0}(40),\chi _{-40}).$
\end{theorem}

\begin{proof}
Appealing to \cite[Theorem 5.9, p.88]{Stein} with $\epsilon =\chi _{-40}\ $%
and $\psi ,$ $\chi \in \{\chi _{1},$ $\chi _{-40},$ $\chi _{-8},$ $\chi
_{5}\}$ we see that $\{E_{3,\chi _{1},\chi _{-40}}(q),E_{3,\chi _{-40},\chi
_{1}}(q),E_{3,\chi _{-8},\chi _{5}}(q),E_{3,\chi _{5},\chi _{-8}}(q)\}$ is a
basis for $E_{3}(\Gamma _{0}(40),\chi _{-40}).$

From (\ref{d1})-(\ref{d10}), (\ref{e1})-(\ref{e3}) and Theorem \ref{thm:1}
we see that all the eta quotients $D_{j}(q)(1\leq j\leq 10)$ are contained
in $S_{3}(\Gamma _{0}(40),\chi _{-40}).$ The eta quotients $D_{j}(q)(1\leq
j\leq 10)$ are linearly independent over $%
\mathbb{C}
$. Since dim($S_{3}(\Gamma _{0}(40),\chi _{-40}))=10$ from (\ref{dim4})$,$ $%
\{D_{j}(q):j=1,2,...,10\}$ form a basis of $S_{3}(\Gamma _{0}(40),\chi
_{-40}).$

Thus, $\{E_{3,\chi _{1},\chi _{-40}}(q),E_{3,\chi _{-40},\chi
_{1}}(q),E_{3,\chi _{-8},\chi _{5}}(q),E_{3,\chi _{5},\chi _{-8}}(q)\}\ $%
together with $\{D_{j}(q):j=1,2,...,10\}$ form a basis of $M_{3}(\Gamma
_{0}(40),\chi _{-40}).$
\end{proof}

\section{Theta Function Identities}

In this section, we give some new theta function identities that will be
used in the next session to prove the main theorem for the representation
number formulae.

\begin{theorem}
\label{thm:6}Let $(l_{1},l_{2},l_{5},l_{10})$ be any of the quadratic forms
listed in the first column of Table \ref{Tab:Tb1}. Then%
\begin{equation*}
\mathop{\displaystyle \prod }\limits_{d\mid 10}\varphi ^{l_{d}}(q^{d})=%
\mathop{\displaystyle \sum }\limits_{t\mid 10}x_{t}E_{3,\chi _{1},\chi
_{-4}}(q^{t})+\mathop{\displaystyle \sum }\limits_{t\mid 10}y_{t}E_{3,\chi
_{-4},\chi _{1}}(q^{t})+\sum_{j=1}^{8}e_{j}A_{j}(q)
\end{equation*}%
where $x_{t},y_{t}$ $(t\mid 10)$ and $e_{j}$ $(1\leq j\leq 8)$ are listed in
Table \ref{Tab:Tb2}.
\end{theorem}

\begin{proof}
It is clear from (\ref{e1}), (\ref{tet3}) and Theorem \ref{thm:1} that $%
\mathop{\displaystyle \prod }\limits_{d\mid 10}\varphi ^{l_{d}}(q^{d})$ is
in $M_{3}(\Gamma _{0}(40),\chi _{-4})$ for each quadruple $%
(l_{1},l_{2},l_{5},l_{10})$ given in the first column of Table \ref{Tab:Tb1}%
. By Theorem \ref{thm:2}, $\mathop{\displaystyle \prod }\limits_{d\mid
10}\varphi ^{l_{d}}(q^{d})$ must be a linear combination of $E_{3,\chi
_{1},\chi _{-4}}(q^{t}),$ $E_{3,\chi _{-4},\chi _{1}}(q^{t})$ $(t\mid 10)$
and $A_{j}(q)(1\leq j\leq 8)$, namely%
\begin{equation}
\mathop{\displaystyle \prod }\limits_{d\mid 10}\varphi ^{l_{d}}(q^{d})=%
\mathop{\displaystyle \sum }\limits_{t\mid 10}x_{t}E_{3,\chi _{1},\chi
_{-4}}(q^{t})+\mathop{\displaystyle \sum }\limits_{t\mid 10}y_{t}E_{3,\chi
_{-4},\chi _{1}}(q^{t})+\sum_{j=1}^{8}e_{j}A_{j}(q)  \label{eqq}
\end{equation}%
for some $x_{t},y_{t}(t\mid 10)$ and $e_{j}(1\leq j\leq 8).$ Appealing to
\cite[Theorem 3.13]{Kilford}, \cite{Sturm} the Sturm bound for $M_{3}(\Gamma
_{0}(40),\chi _{-4})$ is $18$. So, equating the coefficients of $q^{n}$ for $%
1\leq n\leq 18$ on both sides of (\ref{eqq}) we obtain a system of linear
equations with the $16$ unknowns $x_{t},y_{t}(t\mid 10)$ and $e_{j}(1\leq
j\leq 8).$ Solving this system of equations we obtained the desired result$.$
\end{proof}

\begin{theorem}
\label{thm:7}\bigskip Let $(l_{1},l_{2},l_{5},l_{10})$ be any of the
quadratic forms listed in the second column of Table \ref{Tab:Tb1}. Then%
\begin{eqnarray*}
\mathop{\displaystyle \prod }\limits_{d\mid 10}\varphi ^{l_{d}}(q^{d}) &=&%
\mathop{\displaystyle \sum }\limits_{t\mid 2}x_{t}E_{3,\chi _{1},\chi
_{-20}}(q^{t})+\mathop{\displaystyle \sum }\limits_{t\mid 2}y_{t}E_{3,\chi
_{-20},\chi _{1}}(q^{t}) \\
&&+\mathop{\displaystyle \sum }\limits_{t\mid 2}z_{t}E_{3,\chi _{-4},\chi
_{5}}(q^{t})+\mathop{\displaystyle \sum }\limits_{t\mid 2}w_{t}E_{3,\chi
_{5},\chi _{-4}}(q^{t})+\sum_{j=1}^{8}e_{j}B_{j}(q)
\end{eqnarray*}%
where $x_{t},y_{t},z_{t},w_{t}$ $(t\mid 2)$ and $e_{j}$ $(1\leq j\leq 8)$
are listed in Table \ref{Tab:Tb3}.
\end{theorem}

\begin{proof}
It is clear from (\ref{e1}), (\ref{tet3}) and Theorem \ref{thm:1} that $%
\mathop{\displaystyle \prod }\limits_{d\mid 10}\varphi ^{l_{d}}(q^{d})$ is
in $M_{3}(\Gamma _{0}(40),\chi _{-20})$ for each quadruple $%
(l_{1},l_{2},l_{5},l_{10})$ given in the second column of Table \ref{Tab:Tb1}%
. By Theorem \ref{thm:3}, $\mathop{\displaystyle \prod }\limits_{d\mid
10}\varphi ^{l_{d}}(q^{d})$ must be a linear combination of $E_{3,\chi
_{1},\chi _{-20}}(q^{t}),E_{3,\chi _{-20},\chi _{1}}(q^{t}),E_{3,\chi
_{-4},\chi _{5}}(q^{t}),$ $E_{3,\chi _{5},\chi _{-4}}(q^{t})(t\mid 2)\ $and $%
B_{j}(q)(1\leq j\leq 8)$ namely%
\begin{eqnarray}
\mathop{\displaystyle \prod }\limits_{d\mid 10}\varphi ^{l_{d}}(q^{d}) &=&%
\mathop{\displaystyle \sum }\limits_{t\mid 2}x_{t}E_{3,\chi _{1},\chi
_{-20}}(q^{t})+\mathop{\displaystyle \sum }\limits_{t\mid 2}y_{t}E_{3,\chi
_{-20},\chi _{1}}(q^{t})  \label{eqq2} \\
&&+\mathop{\displaystyle \sum }\limits_{t\mid 2}z_{t}E_{3,\chi _{-4},\chi
_{5}}(q^{t})+\mathop{\displaystyle \sum }\limits_{t\mid 2}w_{t}E_{3,\chi
_{5},\chi _{-4}}(q^{t})+\sum_{j=1}^{8}e_{j}B_{j}(q)  \notag
\end{eqnarray}%
for some $x_{t},y_{t},z_{t},w_{t}(t\mid 2)$ and $e_{j}(1\leq j\leq 8).$ So,
equating the coefficients of $q^{n}$ for $1\leq n\leq 18$ on both sides of (%
\ref{eqq2}) we obtain a system of linear equations with the $16$ unknowns.
Solving this system of equations we found values of $%
x_{t},y_{t},z_{t},w_{t}(t\mid 2)$ and $e_{j}(1\leq j\leq 10).$
\end{proof}

\begin{theorem}
\label{thm:8}\bigskip \bigskip Let $(l_{1},l_{2},l_{5},l_{10})$ be any of
the quadratic forms listed in the third column of Table \ref{Tab:Tb1}. Then%
\begin{equation*}
\mathop{\displaystyle \prod }\limits_{d\mid 10}\varphi ^{l_{d}}(q^{d})=%
\mathop{\displaystyle \sum }\limits_{t\mid 5}x_{t}E_{3,\chi _{1},\chi
_{-8}}(q^{t})+\mathop{\displaystyle \sum }\limits_{t\mid 5}y_{t}E_{3,\chi
_{-8},\chi _{1}}(q^{t})+\sum_{j=1}^{10}e_{j}C_{j}(q)
\end{equation*}%
where $x_{t},y_{t}$ $(t=1,5)$ and $e_{j}$ $(1\leq j\leq 10)$ are listed in
Table \ref{Tab:Tb4}.
\end{theorem}

\begin{proof}
It is clear from (\ref{e1}), (\ref{tet3}) and Theorem \ref{thm:1} that $%
\mathop{\displaystyle \prod }\limits_{d\mid 10}\varphi ^{l_{d}}(q^{d})$ is
in $M_{3}(\Gamma _{0}(40),\chi _{-8})$ for each quadruple $%
(l_{1},l_{2},l_{5},l_{10})$ given in the third column of Table \ref{Tab:Tb1}%
. By Theorem \ref{thm:4}, $\mathop{\displaystyle \prod }\limits_{d\mid
10}\varphi ^{l_{d}}(q^{d})$ must be a linear combination of $E_{3,\chi
_{1},\chi _{-8}}(q^{t}),$ $E_{3,\chi _{-8},\chi _{1}}(q^{t})(t\mid 2)\ $and $%
C_{j}(q)(1\leq j\leq 10)$ namely%
\begin{equation}
\mathop{\displaystyle \prod }\limits_{d\mid 10}\varphi ^{l_{d}}(q^{d})=%
\mathop{\displaystyle \sum }\limits_{t\mid 5}x_{t}E_{3,\chi _{1},\chi
_{-8}}(q^{t})+\mathop{\displaystyle \sum }\limits_{t\mid 5}y_{t}E_{3,\chi
_{-8},\chi _{1}}(q^{t})+\sum_{j=1}^{10}e_{j}C_{j}(q)  \label{eqq3}
\end{equation}%
for some $x_{t},y_{t}(t\mid 5)$ and $e_{j}(1\leq j\leq 10).$ Equating the
coefficients of $q^{n}$ for $1\leq n\leq 18$ on both sides of (\ref{eqq3})
we obtain a system of linear equations with the $14$ variables. Solving this
system of equations we found values of $x_{t},y_{t}(t\mid 5)$ and $%
e_{j}(1\leq j\leq 10).$
\end{proof}

\begin{theorem}
\label{thm:9}\bigskip \bigskip Let $(l_{1},l_{2},l_{5},l_{10})$ be any of
the quadratic forms listed in the fourth column of Table \ref{Tab:Tb1}. Then%
\begin{eqnarray*}
\mathop{\displaystyle \prod }\limits_{d\mid 10}\varphi ^{l_{d}}(q^{d})
&=&x_{1}E_{3,\chi _{1},\chi _{-40}}(q)+y_{1}E_{3,\chi _{-40},\chi _{1}}(q) \\
&&+z_{1}E_{3,\chi _{-8},\chi _{5}}(q)+w_{1}E_{3,\chi _{5},\chi
_{-8}}(q)+\sum_{j=1}^{10}e_{j}D_{j}(q)
\end{eqnarray*}%
where $x_{1},y_{1},z_{1},w_{1}$ and $e_{j}$ $(1\leq j\leq 10)$ are listed in
Table \ref{Tab:Tb5}.
\end{theorem}

\begin{proof}
It is clear from (\ref{e1}), (\ref{tet3}) and Theorem \ref{thm:1} that $%
\mathop{\displaystyle \prod }\limits_{d\mid 10}\varphi ^{l_{d}}(q^{d})$ is
in $M_{3}(\Gamma _{0}(40),\chi _{-40})$ for each quadruple $%
(l_{1},l_{2},l_{5},l_{10})$ given in the second column of Table \ref{Tab:Tb1}%
. By Theorem \ref{thm:3}, $\mathop{\displaystyle \prod }\limits_{d\mid
10}\varphi ^{l_{d}}(q^{d})$ must be a linear combination of $E_{3,\chi
_{1},\chi _{-40}}(q),E_{3,\chi _{-40},\chi _{1}}(q),E_{3,\chi _{-8},\chi
_{5}}(q),$ $E_{3,\chi _{5},\chi _{-8}}(q)\ $and $D_{j}(q)(1\leq j\leq 10)$
namely%
\begin{eqnarray}
\mathop{\displaystyle \prod }\limits_{d\mid 10}\varphi ^{l_{d}}(q^{d})
&=&x_{1}E_{3,\chi _{1},\chi _{-40}}(q)+y_{1}E_{3,\chi _{-40},\chi _{1}}(q)
\label{eqq4} \\
&&+z_{1}E_{3,\chi _{-8},\chi _{5}}(q)+w_{1}E_{3,\chi _{5},\chi
_{-8}}(q)+\sum_{j=1}^{10}e_{j}D_{j}(q)  \notag
\end{eqnarray}%
for some $x_{t},y_{t},z_{t},w_{t}$ and $e_{j}(1\leq j\leq 10).$ So, equating
the coefficients of $q^{n}$ for $1\leq n\leq 18$ on both sides of (\ref{eqq4}%
) we obtain a system of linear equations with the $14$ unknowns. Solving
this system of equations we found values of $x_{t},y_{t},z_{t},w_{t}$ and $%
e_{j}(1\leq j\leq 10).$
\end{proof}

\section{Statement of the Main Theorem}

We now give the formulas $N(1^{l_{1}},2^{l_{2}},5^{l_{5}},10^{l_{10}};n)$
for each of the quadratic forms $(l_{1},l_{2},l_{5},l_{10})$ listed in Table
1 in Theorems \ref{thm:10}--\ref{thm:13}.

\begin{theorem}
\label{thm:10}Let $(l_{1},l_{2},l_{5},l_{10})$ be any of the sextenary
quadratic forms listed in the first column of Table \ref{Tab:Tb1}. Then%
\begin{equation*}
N(1^{l_{1}},2^{l_{2}},5^{l_{5}},10^{l_{10}};n)=\mathop{\displaystyle \sum }%
\limits_{t\mid 10}x_{t}\sigma _{(2,\chi _{1},\chi _{-4})}(\frac{n}{t})+%
\mathop{\displaystyle \sum }\limits_{t\mid 10}y_{t}\sigma _{(2,\chi
_{-4},\chi _{1})}(\frac{n}{t})+\sum_{j=1}^{8}e_{j}a_{j}(n)
\end{equation*}%
where $x_{i},y_{i}$ $(t\mid 10)$ and $e_{j}$ $(1\leq j\leq 8)$ are listed in
Table \ref{Tab:Tb2}.
\end{theorem}

\begin{proof}
This follows from (\ref{tet1}), (\ref{Eis2}), (\ref{a9}), and Theorem \ref%
{thm:6}.
\end{proof}

\begin{theorem}
\bigskip \label{thm:11}Let $(l_{1},l_{2},l_{5},l_{10})$ be any of the
quadratic forms listed in the second column of Table \ref{Tab:Tb1}. Then%
\begin{eqnarray*}
N(1^{l_{1}},2^{l_{2}},5^{l_{5}},10^{l_{10}};n) &=&\mathop{\displaystyle \sum
}\limits_{t\mid 2}x_{t}\sigma _{(2,\chi _{1},\chi _{-20})}(\frac{n}{t})+%
\mathop{\displaystyle \sum }\limits_{t\mid 2}y_{t}\sigma _{(2,\chi
_{-20},\chi _{1})}(\frac{n}{t}) \\
&&+\mathop{\displaystyle \sum }\limits_{t\mid 2}z_{t}\sigma _{(2,\chi
_{-4},\chi _{5})}(\frac{n}{t})+\mathop{\displaystyle \sum }\limits_{t\mid
2}w_{t}\sigma _{(2,\chi _{5},\chi _{-4})}(\frac{n}{t})+%
\sum_{j=1}^{8}e_{j}b_{j}(n)
\end{eqnarray*}%
where $x_{i},y_{i}$ $(t\mid 2)$ and $e_{j}$ $(1\leq j\leq 8)$ are listed in
Table \ref{Tab:Tb3}.
\end{theorem}

\begin{proof}
This follows from (\ref{tet1}), (\ref{Eis3})-(\ref{Eis4}), (\ref{b9}), and
Theorem \ref{thm:7}.
\end{proof}

\begin{theorem}
\bigskip \label{thm:12}Let $(l_{1},l_{2},l_{5},l_{10})$ be any of the
quadratic forms listed in the third column of Table \ref{Tab:Tb1}. Then%
\begin{equation*}
N(1^{l_{1}},2^{l_{2}},5^{l_{5}},10^{l_{10}};n)=\mathop{\displaystyle \sum }%
\limits_{t\mid 5}x_{t}\sigma _{(2,\chi _{1},\chi _{-8})}(\frac{n}{t})+%
\mathop{\displaystyle \sum }\limits_{t\mid 5}y_{t}\sigma _{(2,\chi
_{-8},\chi _{1})}(\frac{n}{t})+\sum_{j=1}^{10}e_{j}c_{j}(n)
\end{equation*}%
where $x_{t},y_{t}$ $(t\mid 5)$ and $e_{j}$ $(1\leq j\leq 10)$ are listed in
Table \ref{Tab:Tb4}.
\end{theorem}

\begin{proof}
This follows from (\ref{tet1}), (\ref{Eis5}), (\ref{c11}), and Theorem \ref%
{thm:8}.
\end{proof}

\begin{theorem}
\bigskip \bigskip \label{thm:13}Let $(l_{1},l_{2},l_{5},l_{10})$ be any of
the quadratic forms listed in the fourth column of Table \ref{Tab:Tb1}. Then%
\begin{eqnarray*}
N(1^{l_{1}},2^{l_{2}},5^{l_{5}},10^{l_{10}};n) &=&x_{1}\sigma _{(2,\chi
_{1},\chi _{-40})}(n)+y_{1}\sigma _{(2,\chi _{-40},\chi _{1})}(n) \\
&&+z_{1}\sigma _{(2,\chi _{-8},\chi _{5})}(n)+w_{1}\sigma _{(2,\chi
_{5},\chi _{-8})}(n)+\sum_{j=1}^{10}e_{j}d_{j}(n)
\end{eqnarray*}%
where $x_{1},y_{1},z_{1},w_{1}$ and $e_{j}$ $(1\leq j\leq 10)$ are listed in
Table \ref{Tab:Tb5}.
\end{theorem}

\begin{proof}
This follows from (\ref{tet1}), (\ref{Eis6})-(\ref{Eis7}), (\ref{d11}), and
Theorem \ref{thm:9}.
\end{proof}

\section{Acknowledgement}

This material is based upon work supported by the Simons Foundation
Institute Grant Award ID 507536 while the author was in residence at
the Institute for Computational and Experimental Research in
Mathematics in Providence, RI.

\begin{landscape}

\begin{table}[ht]
\caption{Values of $x_{i},y_{i}$ $(t\mid 10)$ and $e_{j}$ $(1\leq
j\leq 8)$ for Theorem \ref{thm:6} and Theorem \ref{thm:10}.}
\label{Tab:Tb2}\centering

\begin{tabular}{|c|cccccccc|cccccccc|}
\hline
$(l_{1},l_{2},l_{5},l_{10})$ & $x_{1}$ & $x_{2}$ & $x_{5}$ & $x_{10}$ & $y_{1}$ & $y_{2}$ & $%
y_{5}$ & $y_{10}$ & $e_{1}$ & $e_{2}$ & $e_{3}$ & $e_{4}$ & $e_{5}$
& $e_{6}$ & $e_{7}$ & $e_{8}$ \\ \hline\hline
$(6,0,0,0)$ & $-4$ & $0$ & $0$ & $0$ & $16$ & $0$ & $0$ & $0$ & $0$ & $0$ & $%
0$ & $0$ & $0$ & $0$ & $0$ & $0$ \\
$(4,2,0,0)$ & $0$ & $-4$ & $0$ & $0$ & $8$ & $0$ & $0$ & $0$ & $0$ &
$0$ & $0
$ & $0$ & $0$ & $0$ & $0$ & $0$ \\
$(4,0,2,0)$ & $\frac{-24}{31}$ & $0$ & $\frac{-100}{31}$ & $0$ & $\frac{96}{%
31}$ & $0$ & $\frac{400}{31}$ & $0$ & $\frac{-160}{31}$ &
$\frac{624}{31}$ & $\frac{176}{31}$ & $\frac{-1792}{31}$ &
$\frac{1440}{31}$ & $\frac{1344}{31}$
& $\frac{8960}{31}$ & $0$ \\
$(4,0,0,2)$ & $0$ & $\frac{-24}{31}$ & $0$ & $\frac{-100}{31}$ & $\frac{48}{%
31}$ & $0$ & $\frac{200}{31}$ & $0$ & $\frac{-976}{31}$ &
$\frac{1032}{31}$
& $\frac{-24}{31}$ & $\frac{-3840}{31}$ & $\frac{-160}{31}$ & $\frac{3008}{31%
}$ & $\frac{14080}{31}$ & $\frac{9600}{31}$ \\
$(3,1,1,1)$ & $0$ & $\frac{26}{31}$ & $0$ & $\frac{-150}{31}$ & $\frac{52}{31%
}$ & $0$ & $\frac{-300}{31}$ & $0$ & $\frac{-552}{31}$ & $\frac{366}{31}$ & $%
\frac{-78}{31}$ & $\frac{-1696}{31}$ & $\frac{-640}{31}$ &
$\frac{1464}{31}$
& $\frac{5760}{31}$ & $\frac{4480}{31}$ \\
$(2,4,0,0)$ & $0$ & $-4$ & $0$ & $0$ & $4$ & $0$ & $0$ & $0$ & $0$ &
$0$ & $0
$ & $0$ & $0$ & $0$ & $0$ & $0$ \\
$(2,2,2,0)$ & $0$ & $\frac{-24}{31}$ & $0$ & $\frac{-100}{31}$ & $\frac{48}{%
31}$ & $0$ & $\frac{200}{31}$ & $0$ & $\frac{16}{31}$ & $\frac{412}{31}$ & $%
\frac{100}{31}$ & $\frac{-864}{31}$ & $\frac{1080}{31}$ &
$\frac{528}{31}$ &
$\frac{4160}{31}$ & $\frac{-320}{31}$ \\
$(2,2,0,2)$ & $0$ & $\frac{-24}{31}$ & $0$ & $\frac{-100}{31}$ & $\frac{24}{%
31}$ & $0$ & $\frac{100}{31}$ & $0$ & $\frac{-452}{31}$ &
$\frac{556}{31}$ &
$0$ & $\frac{-1840}{31}$ & $\frac{-20}{31}$ & $\frac{1384}{31}$ & $\frac{6560%
}{31}$ & $\frac{4640}{31}$ \\
$(2,0,4,0)$ & $\frac{-4}{31}$ & $0$ & $\frac{-120}{31}$ & $0$ & $\frac{16}{31%
}$ & $0$ & $\frac{480}{31}$ & $0$ & $\frac{32}{155}$ & $\frac{1264}{155}$ & $%
\frac{112}{31}$ & $\frac{-2816}{155}$ & $\frac{736}{31}$ &
$\frac{2112}{155}$
& $\frac{2816}{31}$ & $0$ \\
$(2,0,2,2)$ & $0$ & $\frac{-4}{31}$ & $0$ & $\frac{-120}{31}$ &
$\frac{8}{31} $ & $0$ & $\frac{240}{31}$ & $0$ & $\frac{-896}{155}$
& $\frac{1108}{155}$ & $\frac{-4}{31}$ & $\frac{-4192}{155}$ &
$\frac{56}{31}$ & $\frac{3664}{155}$
& $\frac{3008}{31}$ & $\frac{1600}{31}$ \\
$(2,0,0,4)$ & $0$ & $\frac{-4}{31}$ & $0$ & $\frac{-120}{31}$ &
$\frac{8}{31} $ & $0$ & $\frac{240}{31}$ & $0$ & $\frac{-896}{155}$
& $\frac{1108}{155}$ & $\frac{-4}{31}$ & $\frac{-4192}{155}$ &
$\frac{56}{31}$ & $\frac{3664}{155}$
& $\frac{3008}{31}$ & $\frac{1600}{31}$ \\
$(1,3,1,1)$ & $0$ & $\frac{26}{31}$ & $0$ & $\frac{-150}{31}$ & $\frac{26}{31%
}$ & $0$ & $\frac{-150}{31}$ & $0$ & $\frac{-284}{31}$ & $\frac{88}{31}$ & $%
\frac{-52}{31}$ & $\frac{-604}{31}$ & $\frac{-540}{31}$ &
$\frac{552}{31}$ &
$\frac{2160}{31}$ & $\frac{2000}{31}$ \\
$(1,1,3,1)$ & $0$ & $\frac{6}{31}$ & $0$ & $\frac{-130}{31}$ &
$\frac{12}{31}
$ & $0$ & $\frac{-260}{31}$ & $0$ & $\frac{-656}{155}$ & $\frac{98}{155}$ & $%
\frac{-18}{31}$ & $\frac{-1232}{155}$ & $\frac{-224}{31}$ &
$\frac{1384}{155}
$ & $\frac{1024}{31}$ & $\frac{576}{31}$ \\
$(1,1,1,3)$ & $0$ & $\frac{6}{31}$ & $0$ & $\frac{-130}{31}$ &
$\frac{6}{31}$
& $0$ & $\frac{-130}{31}$ & $0$ & $\frac{-776}{155}$ & $\frac{588}{155}$ & $%
\frac{-12}{31}$ & $\frac{-2452}{155}$ & $\frac{-96}{31}$ &
$\frac{2144}{155}$
& $\frac{1872}{31}$ & $\frac{1072}{31}$ \\
$(0,6,0,0)$ & $0$ & $-4$ & $0$ & $0$ & $0$ & $16$ & $0$ & $0$ & $0$ & $0$ & $%
0$ & $0$ & $0$ & $0$ & $0$ & $0$ \\
$(0,4,2,0)$ & $0$ & $\frac{-24}{31}$ & $0$ & $\frac{-100}{31}$ & $\frac{24}{%
31}$ & $0$ & $\frac{100}{31}$ & $0$ & $\frac{168}{31}$ & $\frac{184}{31}$ & $%
0$ & $\frac{-352}{31}$ & $\frac{600}{31}$ & $\frac{144}{31}$ & $\frac{1600}{%
31}$ & $\frac{-320}{31}$ \\
$(0,4,0,2)$ & $0$ & $\frac{-24}{31}$ & $0$ & $\frac{-100}{31}$ & $0$ & $%
\frac{96}{31}$ & $0$ & $\frac{400}{31}$ & $\frac{-144}{31}$ &
$\frac{256}{31}
$ & $0$ & $\frac{-832}{31}$ & $0$ & $\frac{576}{31}$ & $\frac{2880}{31}$ & $%
\frac{2240}{31}$ \\
$(0,2,4,0)$ & $0$ & $\frac{-4}{31}$ & $0$ & $\frac{-120}{31}$ &
$\frac{8}{31}
$ & $0$ & $\frac{240}{31}$ & $0$ & $\frac{592}{155}$ & $\frac{984}{155}$ & $%
\frac{120}{31}$ & $\frac{-1216}{155}$ & $\frac{800}{31}$ &
$\frac{192}{155}$
& $\frac{1024}{31}$ & $\frac{-384}{31}$ \\
$(0,2,2,2)$ & $0$ & $\frac{-4}{31}$ & $0$ & $\frac{-120}{31}$ &
$\frac{4}{31}
$ & $0$ & $\frac{120}{31}$ & $0$ & $\frac{-108}{155}$ & $\frac{484}{155}$ & $%
0$ & $\frac{-1616}{155}$ & $\frac{100}{31}$ & $\frac{1112}{155}$ & $\frac{928%
}{31}$ & $\frac{608}{31}$ \\
$(0,2,0,4)$ & $0$ & $\frac{-4}{31}$ & $0$ & $\frac{-120}{31}$ &
$\frac{8}{31}
$ & $0$ & $\frac{240}{31}$ & $0$ & $\frac{592}{155}$ & $\frac{984}{155}$ & $%
\frac{120}{31}$ & $\frac{-1216}{155}$ & $\frac{800}{31}$ &
$\frac{192}{155}$
& $\frac{1024}{31}$ & $\frac{-384}{31}$ \\
$(0,0,6,0)$ & $0$ & $0$ & $-4$ & $0$ & $0$ & $0$ & $16$ & $0$ & $0$ & $0$ & $%
0$ & $0$ & $0$ & $0$ & $0$ & $0$ \\
$(0,0,4,2)$ & $0$ & $0$ & $0$ & $-4$ & $0$ & $0$ & $8$ & $0$ & $0$ &
$0$ & $0
$ & $0$ & $0$ & $0$ & $0$ & $0$ \\
$(0,0,2,4)$ & $0$ & $0$ & $0$ & $-4$ & $0$ & $0$ & $4$ & $0$ & $0$ &
$0$ & $0
$ & $0$ & $0$ & $0$ & $0$ & $0$ \\
$(0,0,0,6)$ & $0$ & $0$ & $0$ & $-4$ & $0$ & $0$ & $0$ & $16$ & $0$ & $0$ & $%
0$ & $0$ & $0$ & $0$ & $0$ & $0$ \\ \hline
\end{tabular}
\end{table}

\begin{table}[ht]
\caption{Values of $x_{i},y_{i}$ $(t\mid 2)$ and $e_{j}$ $(1\leq
j\leq 8)$  for Theorem \ref{thm:7} and Theorem \ref{thm:11}.}
\label{Tab:Tb3}\centering

\begin{tabular}{|ccccccccc|cccccccc|}
\hline
\multicolumn{1}{|c|}{$(l_{1},l_{2},l_{5},l_{10})$} & $x_{1}$ & $x_{2}$ & $y_{1}$ & $y_{2}$ & $%
z_{1}$ & $z_{2}$ & $w_{1}$ & $w_{2}$ & $e_{1}$ & $e_{2}$ & $e_{3}$ &
$e_{4}$ & $e_{5}$ & $e_{6}$ & $e_{7}$ & $e_{8}$ \\ \hline\hline
\multicolumn{1}{|c|}{$(5,0,1,0)$} & $\frac{-1}{15}$ & $0$ &
$\frac{20}{3}$ &
$0$ & $\frac{4}{15}$ & $0$ & $\frac{-5}{3}$ & $0$ & $\frac{24}{5}$ & $32$ & $%
0$ & $32$ & $32$ & $0$ & $160$ & $0$ \\
\multicolumn{1}{|c|}{$(4,1,0,1)$} & $0$ & $\frac{-1}{15}$ &
$\frac{10}{3}$ &
$0$ & $\frac{-2}{15}$ & $0$ & $0$ & $\frac{-5}{3}$ & $\frac{-24}{5}$ & $%
\frac{1112}{15}$ & $\frac{48}{5}$ & $\frac{352}{15}$ &
$\frac{-352}{15}$ & $0
$ & $\frac{1120}{3}$ & $\frac{2240}{3}$ \\
\multicolumn{1}{|c|}{$(3,2,1,0)$} & $0$ & $\frac{-1}{15}$ &
$\frac{10}{3}$ &
$0$ & $\frac{2}{15}$ & $0$ & $0$ & $\frac{5}{3}$ & $\frac{34}{15}$ & $\frac{%
232}{15}$ & $\frac{4}{15}$ & $\frac{172}{15}$ & $\frac{308}{15}$ & $\frac{-40%
}{3}$ & $80$ & $0$ \\
\multicolumn{1}{|c|}{$(3,0,3,0)$} & $\frac{-1}{15}$ & $0$ & $\frac{4}{3}$ & $%
0$ & $\frac{4}{15}$ & $0$ & $\frac{-1}{3}$ & $0$ & $\frac{24}{5}$ & $\frac{64%
}{3}$ & $0$ & $\frac{80}{3}$ & $\frac{80}{3}$ & $0$ & $96$ & $0$ \\
\multicolumn{1}{|c|}{$(3,0,1,2)$} & $0$ & $\frac{-1}{15}$ & $\frac{2}{3}$ & $%
0$ & $\frac{2}{15}$ & $0$ & $0$ & $\frac{1}{3}$ & $\frac{-2}{5}$ & $\frac{184%
}{5}$ & $\frac{28}{5}$ & $\frac{124}{5}$ & $\frac{28}{15}$ & $8$ & $\frac{592%
}{3}$ & $\frac{896}{3}$ \\
\multicolumn{1}{|c|}{$(2,3,0,1)$} & $0$ & $\frac{-1}{15}$ & $\frac{5}{3}$ & $%
0$ & $\frac{-1}{15}$ & $0$ & $0$ & $\frac{-5}{3}$ & $\frac{-16}{5}$ & $\frac{%
164}{5}$ & $\frac{28}{5}$ & $\frac{44}{5}$ & $\frac{-212}{15}$ & $0$ & $%
\frac{520}{3}$ & $\frac{1040}{3}$ \\
\multicolumn{1}{|c|}{$(2,1,2,1)$} & $0$ & $\frac{-1}{15}$ & $\frac{2}{3}$ & $%
0$ & $\frac{-2}{15}$ & $0$ & $0$ & $\frac{-1}{3}$ & $\frac{-4}{5}$ & $\frac{%
64}{5}$ & $\frac{64}{15}$ & $\frac{64}{5}$ & $\frac{-32}{15}$ &
$\frac{32}{3}
$ & $\frac{272}{3}$ & $128$ \\
\multicolumn{1}{|c|}{$(2,1,0,3)$} & $0$ & $\frac{-1}{15}$ & $\frac{1}{3}$ & $%
0$ & $\frac{-1}{15}$ & $0$ & $0$ & $\frac{-1}{3}$ & $\frac{-8}{15}$ & $\frac{%
412}{15}$ & $\frac{64}{15}$ & $\frac{232}{15}$ & $\frac{8}{15}$ & $\frac{8}{3%
}$ & $152$ & $208$ \\
\multicolumn{1}{|c|}{$(1,4,1,0)$} & $0$ & $\frac{-1}{15}$ & $\frac{5}{3}$ & $%
0$ & $\frac{1}{15}$ & $0$ & $0$ & $\frac{5}{3}$ & $0$ & $\frac{24}{5}$ & $%
\frac{4}{15}$ & $\frac{4}{5}$ & $\frac{148}{15}$ & $\frac{-40}{3}$ & $\frac{%
80}{3}$ & $0$ \\
\multicolumn{1}{|c|}{$(1,2,3,0)$} & $0$ & $\frac{-1}{15}$ & $\frac{2}{3}$ & $%
0$ & $\frac{2}{15}$ & $0$ & $0$ & $\frac{1}{3}$ & $\frac{18}{5}$ & $\frac{64%
}{5}$ & $\frac{-12}{5}$ & $\frac{44}{5}$ & $\frac{268}{15}$ & $-8$ & $\frac{%
112}{3}$ & $\frac{-64}{3}$ \\
\multicolumn{1}{|c|}{$(1,2,1,2)$} & $0$ & $\frac{-1}{15}$ & $\frac{1}{3}$ & $%
0$ & $\frac{1}{15}$ & $0$ & $0$ & $\frac{1}{3}$ & $0$ & $\frac{272}{15}$ & $%
\frac{8}{5}$ & $\frac{112}{15}$ & $\frac{8}{15}$ & $0$ & $\frac{256}{3}$ & $%
\frac{416}{3}$ \\
\multicolumn{1}{|c|}{$(1,0,5,0)$} & $\frac{-1}{15}$ & $0$ &
$\frac{4}{15}$ &
$0$ & $\frac{4}{15}$ & $0$ & $\frac{-1}{15}$ & $0$ & $\frac{8}{5}$ & $\frac{%
32}{5}$ & $0$ & $\frac{32}{5}$ & $\frac{32}{5}$ & $0$ &
$\frac{96}{5}$ & $0$
\\
\multicolumn{1}{|c|}{$(1,0,3,2)$} & $0$ & $\frac{-1}{15}$ &
$\frac{2}{15}$ &
$0$ & $\frac{2}{15}$ & $0$ & $0$ & $\frac{1}{15}$ & $\frac{-2}{15}$ & $\frac{%
64}{15}$ & $\frac{28}{15}$ & $\frac{76}{15}$ & $\frac{4}{3}$ &
$\frac{8}{3}$
& $\frac{144}{5}$ & $\frac{192}{5}$ \\
\multicolumn{1}{|c|}{$(1,0,1,4)$} & $0$ & $\frac{-1}{15}$ &
$\frac{1}{15}$ &
$0$ & $\frac{1}{15}$ & $0$ & $0$ & $\frac{1}{15}$ & $0$ & $\frac{24}{5}$ & $%
\frac{28}{15}$ & $\frac{28}{5}$ & $\frac{28}{15}$ & $\frac{8}{3}$ & $\frac{%
496}{15}$ & $\frac{192}{5}$ \\
\multicolumn{1}{|c|}{$(0,5,0,1)$} & $0$ & $\frac{-1}{15}$ & $0$ & $\frac{20}{%
3}$ & $0$ & $\frac{4}{15}$ & $0$ & $\frac{-5}{3}$ & $\frac{-16}{5}$ & $\frac{%
48}{5}$ & $\frac{16}{5}$ & $\frac{8}{5}$ & $\frac{-48}{5}$ & $0$ & $80$ & $%
160$ \\
\multicolumn{1}{|c|}{$(0,3,2,1)$} & $0$ & $\frac{-1}{15}$ & $\frac{1}{3}$ & $%
0$ & $\frac{-1}{15}$ & $0$ & $0$ & $\frac{-1}{3}$ & $\frac{-8}{15}$ & $\frac{%
52}{15}$ & $\frac{4}{15}$ & $\frac{52}{15}$ & $\frac{-52}{15}$ &
$\frac{8}{3}
$ & $32$ & $48$ \\
\multicolumn{1}{|c|}{$(0,3,0,3)$} & $0$ & $\frac{-1}{15}$ & $0$ & $\frac{4}{3%
}$ & $0$ & $\frac{4}{15}$ & $0$ & $\frac{-1}{3}$ & $\frac{-8}{15}$ & $\frac{%
184}{15}$ & $\frac{8}{15}$ & $\frac{64}{15}$ & $\frac{-8}{5}$ &
$\frac{-8}{3}
$ & $\frac{200}{3}$ & $96$ \\
\multicolumn{1}{|c|}{$(0,1,4,1)$} & $0$ & $\frac{-1}{15}$ &
$\frac{2}{15}$ & $0$ & $\frac{-2}{15}$ & $0$ & $0$ & $\frac{-1}{15}$
& $0$ & $\frac{8}{15}$ &
$0$ & $\frac{16}{15}$ & $\frac{-16}{15}$ & $0$ & $\frac{32}{15}$ & $\frac{64%
}{15}$ \\
\multicolumn{1}{|c|}{$(0,1,2,3)$} & $0$ & $\frac{-1}{15}$ &
$\frac{1}{15}$ & $0$ & $\frac{-1}{15}$ & $0$ & $0$ & $\frac{-1}{15}$
& $0$ & $\frac{12}{5}$ &
$0$ & $\frac{8}{5}$ & $\frac{-8}{15}$ & $0$ & $\frac{176}{15}$ & $\frac{304}{%
15}$ \\
\multicolumn{1}{|c|}{$(0,1,0,5)$} & $0$ & $\frac{-1}{15}$ & $0$ & $\frac{4}{%
15}$ & $0$ & $\frac{4}{15}$ & $0$ & $\frac{-1}{15}$ & $0$ &
$\frac{16}{5}$ & $0$ & $\frac{8}{5}$ & $0$ & $0$ & $16$ &
$\frac{96}{5}$ \\ \hline
\end{tabular}
\end{table}

\begin{table}[ht]
\caption{Values of $x_{t},y_{t}$ $(t\mid 5)$ and $e_{r}$ $(1\leq
j\leq 10)$ for Theorem \ref{thm:8} and Theorem \ref{thm:12}.}
\label{Tab:Tb4}\centering

\begin{tabular}{|c|cccc|cccccccccc|}
\hline
$(l_{1},l_{2},l_{5},l_{10})$ & $x_{1}$ & $x_{5}$ & $y_{1}$ & $y_{5}$ & $e_{1}$ & $e_{2}$ & $%
e_{3}$ & $e_{4}$ & $e_{5}$ & $e_{6}$ & $e_{7}$ & $e_{8}$ & $e_{9}$ &
$e_{10}$
\\ \hline\hline
$(5,1,0,0)$ & $\frac{-2}{3}$ & $0$ & $\frac{32}{3}$ & $0$ & $0$ &
$0$ & $0$
& $0$ & $0$ & $0$ & $0$ & $0$ & $0$ & $0$ \\
$(4,0,1,1)$ & $\frac{8}{63}$ & $\frac{-50}{63}$ & $\frac{128}{63}$ & $\frac{%
800}{63}$ & $\frac{-2512}{63}$ & $\frac{1136}{63}$ & $\frac{-2176}{21}$ & $%
\frac{-3008}{63}$ & $\frac{5552}{63}$ & $\frac{3376}{63}$ &
$\frac{2560}{21}$
& $\frac{2560}{7}$ & $\frac{2560}{21}$ & $\frac{320}{21}$ \\
$(3,3,0,0)$ & $\frac{-2}{3}$ & $0$ & $\frac{16}{3}$ & $0$ &
$\frac{-8}{3}$ & $\frac{4}{3}$ & $-12$ & $\frac{-4}{3}$ &
$\frac{4}{3}$ & $\frac{8}{3}$ & $0$
& $0$ & $0$ & $0$ \\
$(3,1,2,0)$ & $\frac{-26}{189}$ & $\frac{-100}{189}$ & $\frac{416}{189}$ & $%
\frac{-1600}{189}$ & $\frac{-4744}{189}$ & $\frac{3400}{189}$ & $\frac{-352}{%
27}$ & $\frac{-2816}{63}$ & $\frac{392}{3}$ & $\frac{3064}{63}$ & $\frac{%
11264}{63}$ & $\frac{102080}{189}$ & $\frac{55040}{189}$ &
$\frac{-160}{189}$
\\
$(3,1,0,2)$ & $\frac{-26}{189}$ & $\frac{-100}{189}$ & $\frac{208}{189}$ & $%
\frac{-800}{189}$ & $\frac{-3344}{189}$ & $\frac{1816}{189}$ & $\frac{-764}{%
21}$ & $\frac{-556}{27}$ & $\frac{9664}{189}$ & $\frac{692}{27}$ & $\frac{%
4352}{63}$ & $\frac{14240}{63}$ & $\frac{8320}{63}$ & $\frac{320}{21}$ \\
$(2,2,1,1)$ & $\frac{8}{63}$ & $\frac{-50}{63}$ & $\frac{64}{63}$ & $\frac{%
400}{63}$ & $\frac{-776}{63}$ & $\frac{332}{63}$ & $\frac{-1676}{63}$ & $%
\frac{-88}{7}$ & $\frac{548}{21}$ & $\frac{108}{7}$ & $\frac{640}{21}$ & $%
\frac{7360}{63}$ & $\frac{2560}{63}$ & $\frac{64}{9}$ \\
$(2,0,3,1)$ & $\frac{4}{189}$ & $\frac{-130}{189}$ & $\frac{64}{189}$ & $%
\frac{2080}{189}$ & $\frac{-22408}{945}$ & $\frac{2080}{189}$ & $\frac{-20896%
}{315}$ & $\frac{-25664}{945}$ & $\frac{45632}{945}$ &
$\frac{29104}{945}$ & $\frac{22528}{315}$ & $\frac{3968}{21}$ &
$\frac{3968}{63}$ & $\frac{160}{63}
$ \\
$(2,0,1,3)$ & $\frac{4}{189}$ & $\frac{-130}{189}$ & $\frac{32}{189}$ & $%
\frac{1040}{189}$ & $\frac{-20168}{945}$ & $\frac{10112}{945}$ & $\frac{%
-52532}{945}$ & $\frac{-2656}{105}$ & $\frac{15776}{315}$ &
$\frac{3056}{105}
$ & $\frac{23104}{315}$ & $\frac{37280}{189}$ & $\frac{14720}{189}$ & $\frac{%
544}{135}$ \\
$(1,5,0,0)$ & $\frac{-2}{3}$ & $0$ & $\frac{8}{3}$ & $0$ & $0$ & $0$
& $0$ &
$0$ & $0$ & $0$ & $0$ & $0$ & $0$ & $0$ \\
$(1,3,2,0)$ & $\frac{-26}{189}$ & $\frac{-100}{189}$ & $\frac{208}{189}$ & $%
\frac{-800}{189}$ & $\frac{-320}{189}$ & $\frac{1060}{189}$ &
$\frac{580}{21}
$ & $\frac{-340}{27}$ & $\frac{10420}{189}$ & $\frac{368}{27}$ & $\frac{4352%
}{63}$ & $\frac{14240}{63}$ & $\frac{8320}{63}$ & $\frac{-16}{21}$ \\
$(1,3,0,2)$ & $\frac{-26}{189}$ & $\frac{-100}{189}$ & $\frac{104}{189}$ & $%
\frac{-400}{189}$ & $\frac{-376}{189}$ & $\frac{268}{189}$ &
$\frac{368}{189} $ & $\frac{-160}{63}$ & $\frac{716}{63}$ &
$\frac{260}{63}$ & $\frac{128}{9}$
& $\frac{13040}{189}$ & $\frac{9920}{189}$ & $\frac{1376}{189}$ \\
$(1,1,4,0)$ & $\frac{-2}{63}$ & $\frac{-40}{63}$ & $\frac{32}{63}$ & $\frac{%
-640}{63}$ & $\frac{-688}{63}$ & $\frac{352}{63}$ & $\frac{-256}{9}$ & $%
\frac{-320}{21}$ & $32$ & $\frac{352}{21}$ & $\frac{1280}{21}$ & $\frac{8576%
}{63}$ & $\frac{5888}{63}$ & $\frac{-64}{63}$ \\
$(1,1,2,2)$ & $\frac{-2}{63}$ & $\frac{-40}{63}$ & $\frac{16}{63}$ & $\frac{%
-320}{63}$ & $\frac{-2488}{315}$ & $\frac{256}{63}$ & $\frac{-772}{35}$ & $%
\frac{-452}{45}$ & $\frac{6392}{315}$ & $\frac{532}{45}$ &
$\frac{3328}{105}$
& $\frac{1664}{21}$ & $\frac{640}{21}$ & $\frac{16}{7}$ \\
$(1,1,0,4)$ & $\frac{-2}{63}$ & $\frac{-40}{63}$ & $\frac{8}{63}$ & $\frac{%
-160}{63}$ & $\frac{-2768}{315}$ & $\frac{1544}{315}$ &
$\frac{-3296}{315}$
& $\frac{-1328}{105}$ & $\frac{3448}{105}$ & $\frac{1528}{105}$ & $\frac{128%
}{3}$ & $\frac{9248}{63}$ & $\frac{3968}{63}$ & $\frac{736}{315}$ \\
$(0,4,1,1)$ & $\frac{8}{63}$ & $\frac{-50}{63}$ & $\frac{32}{63}$ & $\frac{%
200}{63}$ & $\frac{-160}{63}$ & $\frac{8}{9}$ & $\frac{-128}{21}$ & $\frac{64%
}{63}$ & $\frac{-184}{63}$ & $\frac{-104}{63}$ & $\frac{-320}{21}$ & $\frac{%
-160}{21}$ & $0$ & $\frac{64}{21}$ \\
$(0,2,3,1)$ & $\frac{4}{189}$ & $\frac{-130}{189}$ & $\frac{32}{189}$ & $%
\frac{1040}{189}$ & $\frac{-3536}{945}$ & $\frac{1796}{945}$ & $\frac{-4148}{%
945}$ & $\frac{-472}{105}$ & $\frac{2924}{315}$ & $\frac{452}{105}$ & $\frac{%
2944}{315}$ & $\frac{7040}{189}$ & $\frac{2624}{189}$ & $\frac{112}{135}$ \\
$(0,2,1,3)$ & $\frac{4}{189}$ & $\frac{-130}{189}$ & $\frac{16}{189}$ & $%
\frac{520}{189}$ & $\frac{-2416}{945}$ & $\frac{236}{135}$ &
$\frac{-64}{63}$
& $\frac{-3368}{945}$ & $\frac{1924}{189}$ & $\frac{3268}{945}$ & $\frac{3232%
}{315}$ & $\frac{2608}{63}$ & $\frac{64}{3}$ & $\frac{496}{315}$ \\
$(0,0,5,1)$ & $0$ & $\frac{-2}{3}$ & $0$ & $\frac{32}{3}$ & $0$ &
$0$ & $0$
& $0$ & $0$ & $0$ & $0$ & $0$ & $0$ & $0$ \\
$(0,0,3,3)$ & $0$ & $\frac{-2}{3}$ & $0$ & $\frac{16}{3}$ &
$\frac{16}{75}$ & $\frac{-16}{75}$ & $\frac{-92}{75}$ &
$\frac{16}{25}$ & $\frac{-64}{25}$ &
$\frac{-16}{25}$ & $\frac{-64}{25}$ & $\frac{-32}{3}$ & $\frac{-64}{15}$ & $%
\frac{16}{75}$ \\
$(0,0,1,5)$ & $0$ & $\frac{-2}{3}$ & $0$ & $\frac{8}{3}$ & $0$ & $0$
& $0$ & $0$ & $0$ & $0$ & $0$ & $0$ & $0$ & $0$ \\ \hline
\end{tabular}
\end{table}

\begin{table}[ht]
 \caption{Values of $x_{1},y_{1},z_{1},w_{1}$ and $e_{j}$
$(1\leq j\leq 10)$ for Theorem \ref{thm:9} and Theorem \ref{thm:13}}
\label{Tab:Tb5} \centering

\begin{tabular}{|c|cccc|cccccccccc|}
\hline
$(l_{1},l_{2},l_{5},l_{10})$ & $x_{1}$ & $y_{1}$ & $z_{1}$ & $t_{1}$ & $e_{1}$ & $e_{2}$ & $%
e_{3}$ & $e_{4}$ & $e_{5}$ & $e_{6}$ & $e_{7}$ & $e_{8}$ & $e_{9}$ &
$e_{10}$
\\ \hline\hline
\multicolumn{1}{|c|}{$(5,0,0,1)$} & $\frac{-1}{79}$ & $\frac{400}{79}$ & $%
\frac{16}{79}$ & $\frac{-25}{79}$ & $\frac{1840}{79}$ & $\frac{-800}{79}$ & $%
\frac{320}{79}$ & $0$ & $\frac{-160}{79}$ & $\frac{-1680}{79}$ & $\frac{400}{%
79}$ & $\frac{3040}{79}$ & $0$ & $\frac{7600}{79}$ \\
\multicolumn{1}{|c|}{$(4,1,1,0)$} & $\frac{-1}{79}$ & $\frac{400}{79}$ & $%
\frac{-16}{79}$ & $\frac{25}{79}$ & $\frac{-6112}{79}$ &
$\frac{3136}{79}$ & $\frac{-5312}{79}$ & $\frac{640}{79}$ &
$\frac{2560}{79}$ & $\frac{5280}{79}$
& $\frac{224}{79}$ & $\frac{480}{79}$ & $\frac{12800}{79}$ & $0$ \\
\multicolumn{1}{|c|}{$(3,2,0,1)$} & $\frac{-1}{79}$ & $\frac{200}{79}$ & $%
\frac{8}{79}$ & $\frac{-25}{79}$ & $\frac{1732}{79}$ & $\frac{-916}{79}$ & $%
\frac{1012}{79}$ & $0$ & $\frac{-480}{79}$ & $\frac{-1940}{79}$ & $\frac{292%
}{79}$ & $\frac{1440}{79}$ & $0$ & $\frac{3600}{79}$ \\
\multicolumn{1}{|c|}{$(3,0,2,1)$} & $\frac{-1}{79}$ & $\frac{80}{79}$ & $%
\frac{16}{79}$ & $\frac{-5}{79}$ & $\frac{-408}{79}$ & $\frac{224}{79}$ & $%
\frac{-1792}{79}$ & $\frac{-248}{79}$ & $\frac{896}{79}$ &
$\frac{752}{79}$ & $\frac{384}{79}$ & $\frac{2128}{79}$ &
$\frac{4960}{79}$ & $\frac{1400}{79}
$ \\
\multicolumn{1}{|c|}{$(3,0,0,3)$} & $\frac{-1}{79}$ & $\frac{40}{79}$ & $%
\frac{8}{79}$ & $\frac{-5}{79}$ & $\frac{-1812}{79}$ & $\frac{924}{79}$ & $%
\frac{-2676}{79}$ & $\frac{36}{79}$ & $\frac{1344}{79}$ &
$\frac{1632}{79}$ & $\frac{432}{79}$ & $\frac{2088}{79}$ &
$\frac{8160}{79}$ & $\frac{2160}{79}
$ \\
\multicolumn{1}{|c|}{$(2,3,1,0)$} & $\frac{-1}{79}$ & $\frac{200}{79}$ & $%
\frac{-8}{79}$ & $\frac{25}{79}$ & $\frac{-2764}{79}$ & $\frac{1528}{79}$ & $%
\frac{-2064}{79}$ & $\frac{720}{79}$ & $\frac{960}{79}$ &
$\frac{1900}{79}$
& $\frac{100}{79}$ & $\frac{-240}{79}$ & $\frac{4800}{79}$ & $0$ \\
\multicolumn{1}{|c|}{$(2,1,3,0)$} & $\frac{-1}{79}$ & $\frac{80}{79}$ & $%
\frac{-16}{79}$ & $\frac{5}{79}$ & $\frac{-5000}{79}$ & $\frac{2640}{79}$ & $%
\frac{-4480}{79}$ & $\frac{384}{79}$ & $\frac{2224}{79}$ &
$\frac{3976}{79}$ & $\frac{248}{79}$ & $\frac{208}{79}$ &
$\frac{7680}{79}$ & $\frac{-200}{79}$
\\
\multicolumn{1}{|c|}{$(2,1,1,2)$} & $\frac{-1}{79}$ & $\frac{40}{79}$ & $%
\frac{-8}{79}$ & $\frac{5}{79}$ & $\frac{-1052}{79}$ & $\frac{584}{79}$ & $%
\frac{-1544}{79}$ & $\frac{116}{79}$ & $\frac{760}{79}$ &
$\frac{472}{79}$ & $\frac{280}{79}$ & $\frac{1048}{79}$ &
$\frac{2880}{79}$ & $\frac{1400}{79}$
\\
\multicolumn{1}{|c|}{$(1,4,0,1)$} & $\frac{-1}{79}$ & $\frac{100}{79}$ & $%
\frac{4}{79}$ & $\frac{-25}{79}$ & $\frac{1520}{79}$ & $\frac{-816}{79}$ & $%
\frac{1200}{79}$ & $0$ & $\frac{-640}{79}$ & $\frac{-1280}{79}$ & $\frac{80}{%
79}$ & $\frac{640}{79}$ & $0$ & $\frac{1600}{79}$ \\
\multicolumn{1}{|c|}{$(1,2,2,1)$} & $\frac{-1}{79}$ & $\frac{40}{79}$ & $%
\frac{8}{79}$ & $\frac{-5}{79}$ & $\frac{400}{79}$ & $\frac{-340}{79}$ & $%
\frac{-148}{79}$ & $\frac{-280}{79}$ & $\frac{80}{79}$ & $\frac{52}{79}$ & $%
\frac{116}{79}$ & $\frac{824}{79}$ & $\frac{1840}{79}$ & $\frac{580}{79}$ \\
\multicolumn{1}{|c|}{$(1,2,0,3)$} & $\frac{-1}{79}$ & $\frac{20}{79}$ & $%
\frac{4}{79}$ & $\frac{-5}{79}$ & $\frac{-460}{79}$ & $\frac{168}{79}$ & $%
\frac{-748}{79}$ & $\frac{20}{79}$ & $\frac{304}{79}$ & $\frac{492}{79}$ & $%
\frac{140}{79}$ & $\frac{804}{79}$ & $\frac{3440}{79}$ & $\frac{960}{79}$ \\
\multicolumn{1}{|c|}{$(1,0,4,1)$} & $\frac{-1}{79}$ & $\frac{16}{79}$ & $%
\frac{16}{79}$ & $\frac{-1}{79}$ & $\frac{768}{395}$ & $\frac{-384}{395}$ & $%
\frac{-192}{79}$ & $\frac{-224}{395}$ & $\frac{96}{79}$ &
$\frac{-128}{395}$ & $\frac{128}{79}$ & $\frac{2144}{395}$ &
$\frac{896}{79}$ & $\frac{160}{79}$
\\
\multicolumn{1}{|c|}{$(1,0,2,3)$} & $\frac{-1}{79}$ & $\frac{8}{79}$ & $%
\frac{8}{79}$ & $\frac{-1}{79}$ & $\frac{-56}{79}$ & $\frac{28}{79}$ & $%
\frac{-380}{79}$ & $\frac{-20}{79}$ & $\frac{192}{79}$ & $\frac{8}{79}$ & $%
\frac{144}{79}$ & $\frac{448}{79}$ & $\frac{944}{79}$ & $\frac{292}{79}$ \\
\multicolumn{1}{|c|}{$(1,0,0,5)$} & $\frac{-1}{79}$ & $\frac{4}{79}$ & $%
\frac{4}{79}$ & $\frac{-1}{79}$ & $\frac{-224}{79}$ &
$\frac{112}{79}$ & $-8$
& $\frac{24}{79}$ & $\frac{240}{79}$ & $\frac{88}{79}$ & $\frac{152}{79}$ & $%
\frac{584}{79}$ & $\frac{1600}{79}$ & $\frac{200}{79}$ \\
\multicolumn{1}{|c|}{$(0,5,1,0)$} & $\frac{-1}{79}$ & $\frac{100}{79}$ & $%
\frac{-4}{79}$ & $\frac{25}{79}$ & $\frac{-1880}{79}$ & $\frac{1040}{79}$ & $%
\frac{-440}{79}$ & $\frac{760}{79}$ & $\frac{160}{79}$ &
$\frac{1000}{79}$ &
$\frac{-120}{79}$ & $\frac{-600}{79}$ & $\frac{800}{79}$ & $0$ \\
\multicolumn{1}{|c|}{$(0,3,3,0)$} & $\frac{-1}{79}$ & $\frac{40}{79}$ & $%
\frac{-8}{79}$ & $\frac{5}{79}$ & $\frac{-3264}{79}$ & $\frac{1848}{79}$ & $%
\frac{-2808}{79}$ & $\frac{432}{79}$ & $\frac{1392}{79}$ &
$\frac{2052}{79}$ & $\frac{-36}{79}$ & $\frac{-216}{79}$ &
$\frac{2880}{79}$ & $\frac{-180}{79}
$ \\
\multicolumn{1}{|c|}{$(0,3,1,2)$} & $\frac{-1}{79}$ & $\frac{20}{79}$ & $%
\frac{-4}{79}$ & $\frac{5}{79}$ & $\frac{-816}{79}$ & $\frac{504}{79}$ & $%
\frac{-708}{79}$ & $\frac{140}{79}$ & $\frac{344}{79}$ & $\frac{300}{79}$ & $%
\frac{-20}{79}$ & $\frac{204}{79}$ & $\frac{480}{79}$ & $\frac{620}{79}$ \\
\multicolumn{1}{|c|}{$(0,1,5,0)$} & $\frac{-1}{79}$ & $\frac{16}{79}$ & $%
\frac{-16}{79}$ & $\frac{1}{79}$ & $\frac{-1744}{79}$ & $\frac{1024}{79}$ & $%
\frac{-1280}{79}$ & $\frac{80}{79}$ & $\frac{640}{79}$ &
$\frac{1440}{79}$ &
$0$ & $\frac{-352}{79}$ & $\frac{1600}{79}$ & $\frac{-240}{79}$ \\
\multicolumn{1}{|c|}{$(0,1,3,2)$} & $\frac{-1}{79}$ & $\frac{8}{79}$ & $%
\frac{-8}{79}$ & $\frac{1}{79}$ & $\frac{-704}{395}$ & $\frac{712}{395}$ & $%
\frac{-176}{79}$ & $\frac{292}{395}$ & $\frac{88}{79}$ &
$\frac{-16}{395}$ &
$0$ & $\frac{208}{395}$ & $\frac{-32}{79}$ & $\frac{100}{79}$ \\
\multicolumn{1}{|c|}{$(0,1,1,4)$} & $\frac{-1}{79}$ & $\frac{4}{79}$ & $%
\frac{-4}{79}$ & $\frac{1}{79}$ & $\frac{-224}{79}$ & $\frac{144}{79}$ & $%
\frac{-256}{79}$ & $\frac{16}{79}$ & $\frac{128}{79}$ &
$\frac{160}{79}$ & $0 $ & $\frac{112}{79}$ & $\frac{416}{79}$ &
$\frac{112}{79}$ \\ \hline
\end{tabular}

\end{table}

\end{landscape}

\end{document}